\documentclass[12pt]{article}
\usepackage{latexsym}
\newtheorem{theorem}{Theorem}
\newtheorem{cor}{Corollary}
\title{BG-ranks and 2-cores}
\author{William Y. C. Chen, Kathy Q. Ji, and Herbert S. Wilf}
\begin{document}
\maketitle
\begin{abstract}
We find the number of partitions of $n$ whose BG-rank is $j$, in terms of $pp(n)$, the number of pairs of partitions whose total number of cells is $n$, giving both bijective and generating function proofs. Next we find congruences mod 5 for $pp(n)$, and then we use these to give a new proof of a refined system of congruences for $p(n)$ that was found by Berkovich and Garvan.
\end{abstract}
\section{Introduction} If $\pi$ is a partition of $n$ we define the
\textit{BG-rank} $\beta(\pi)$, of $\pi$ as follows. First draw the
Ferrers diagram of $\pi$. Then fill the cells with alternating
$\pm 1$'s, chessboard style, beginning with a $+1$ in the $(1,1)$
position. The sum of these entries is $\beta(\pi)$, the BG-rank of
$\pi$. For example, the BG-rank of the partition $13=4+3+3+1+1+1$
is $-1$.

\begin{figure}
\centering
\setlength{\unitlength}{.3in}
\begin{picture}(8,6)(-2,0)
\put(.2,5.3){+1}\put(1.2,5.3){$-1$}\put(2.2,5.3){+1}\put(3.2,5.3){$-1$}
\put(.2,4.3){$-1$}\put(1.2,4.3){+1}\put(2.2,4.3){$-1$}
\put(.2,3.3){+1}\put(1.2,3.3){$-1$}\put(2.2,3.3){+1}
\put(.2,2.3){$-1$}
\put(.2,1.3){+1}
\put(.2,0.3){$-1$}
\put(0,6){\line(1,0){4}}
\put(0,5){\line(1,0){4}}
\put(0,4){\line(1,0){3}}
\put(0,3){\line(1,0){3}}
\put(0,2){\line(1,0){1}}
\put(0,1){\line(1,0){1}}
\put(0,0){\line(1,0){1}}
\put(0,6){\line(0,-1){6}}
\put(1,6){\line(0,-1){6}}
\put(2,6){\line(0,-1){3}}
\put(3,6){\line(0,-1){3}}
\put(4,6){\line(0,-1){1}}
\end{picture}
\caption{A partition with BG-rank $-1$}
\end{figure}

This partition statistic has been encountered by several authors
(\cite{bg, bg2, dnm,sch,syd}), but its systematic study was
initiated in \cite{bg}. Here we wish to study the function
\[p_j(n)=\left\vert\, \{\pi:\,\vert\pi\vert= n\,\mathrm{and}\,\beta(\pi)=j\,\}\,\right\vert .\]
We will find a fairly explicit formula for it (see (\ref{eq:fst}) below), and a bijective proof for this formula. We will then show that a number of congruences from \cite{bg} can all be proved from a single set of congruences for the function $pp(n)$ defined by (\ref{eq:pp}) below.
\section{The theorem}
We write $p(n)$ for the usual partition function, and ${\cal
P}(x)$ for its generating function. If $\pi$ is a partition of $n$
then we will write $\vert \pi\vert=n$. $pp(n)$ will be the number
of ordered pairs $\pi',\pi''$ of partitions such that $\vert
\pi'\vert+\vert\pi''\vert=n$, i.e., $pp(n)$ is the sequence that
is generated by
\begin{equation}
\label{eq:pp}
\sum_{n\ge 0}pp(n)x^n={\cal P}(x)^2=\prod_{i\ge 1}\frac{1}{(1-x^i)^2}.
\end{equation}
By convention $pp(n)$ vanishes
unless its argument is a nonnegative integer. Our main result is
as follows.

\begin{theorem} The number of partitions of $n$ whose BG-rank is
$j$ is given by
\begin{equation}\label{eq:fst}p_j(n)=pp\left(\frac{n-j(2j-1)}{2}\right).\end{equation}
\end{theorem}

A non-bijective proof of this is easy, given the results of
\cite{bg}. The authors of \cite{bg} found the two variable
generating function for $\bar p_j(m,n)$, the number of partitions
of $n$ with BG-rank $=j$ and ``2-quotient-rank'' $=m$, in the form
\[\sum_{n,m}\bar
p_j(m,n)x^mq^n=\frac{q^{j(2j-1)}}{(q^2x,q^2/x;q^2)_{\infty}}.\] If
we simply put $x=1$ here, and read off the coefficients of like
powers of $q$, we have (\ref{eq:fst}). $\Box$

\section{Bijective proof}
 A bijective proof of (\ref{eq:fst}) follows from the theory of 2-cores. The \textit{2-core} of a partition $\pi$ is obtained as follows.
Begin with the Ferrers diagram of $\pi$. Then delete a horizontal
or a vertical pair of adjacent cells, subject only to the
restriction that the result of the deletion must be a valid
Ferrers diagram. Repeat this process, making arbitrary choices,
until no further such deletions are possible. The remaining
diagram is the 2-core of $\pi$, $C(\pi)$, say.

The 2-core of a partition is always a staircase partition, i.e., a
partition of the form
 \[{k+1\choose 2}=k+(k-1)+\dots+1.\]
The following representation theorem is well known, and probably
goes back to Littlewood \cite{lit} or to Nakayama \cite{nak}. For
a lucid exposition see Schmidt \cite{sch}.

\begin{theorem}
There is a 1-1 (constructive) correspondence between partitions
$\pi$ of $n$ and triples $(S,\pi',\pi'')$, where $S$ is a
staircase partition (the 2-core of $\pi$), and $\pi',\pi''$ are
partitions such that $n=\vert
S\vert+2\vert\pi'\vert+2\vert\pi''\vert$.
\end{theorem}

The proof of Theorem 1 will follow from the following
observations:
\begin{enumerate} \item First, the BG-rank of a partition and of
its 2-core are equal, since at each stage of the construction of
the 2-core we delete a pair of adjacent cells, which does not
change the BG-rank. \item An easy calculation shows that the
BG-rank of a staircase partition of height $k$ is $(k+1)/2$, if
$k$ is odd, and $-k/2$, if $k$ is even. \item Therefore, if $\pi$
is a partition of BG-rank $=j$ then its 2-core is a staircase
partition of height $2j-1$, if $j>0$, and $-2j$, if $j\le 0$.
\item In either case, if $\pi$ is a partition whose BG-rank is
$j$, then its 2-core is a diagram of exactly $j(2j-1)$ cells,
i.e., a partition of the integer $j(2j-1)$.
\end{enumerate}
Theorem 1 now follows from Theorem 2 and remark 4 above. $\Box$

\begin{cor}
There exists a partition of $n$ with BG-rank $=j$ if and only if
$j+n$ is even and $j(2j-1)\le n$.\end{cor}
\section{Congruences}
The motivation for introducing the BG-rank lay in the wish to refine some known congruences for $p(n)$. We can give quite elementary proofs of some of their congruences, in particular the following:
\begin{eqnarray}\label{eq:p0}p_j(5n)&\equiv& 0\ (\mathrm{mod}\, 5), \ \mathrm{ if}\  j\equiv 1,2 \,(\mathrm{mod}\, 5),\\
\label{eq:p1}p_j(5n+1)&\equiv& 0\ (\mathrm{mod} \, 5),\ \mathrm{ if}\ j\equiv 0,3,4 \,(\mathrm{mod} \, 5),\\
\label{eq:p2}p_j(5n+2)&\equiv& 0\ (\mathrm{mod} \, 5), \ \mathrm{ if}\ j\equiv 1,2,4 \,(\mathrm{mod} \, 5),\\
\label{eq:p3}p_j(5n+3)&\equiv& 0\ (\mathrm{mod} \, 5), \ \mathrm{ if}\  j\equiv 0,3 \,(\mathrm{mod} \, 5),\\
\label{eq:p4}p_j(5n+4)&\equiv& 0\ (\mathrm{mod}\, 5),\ \forall\, j.
\end{eqnarray}

First, we claim that all of the above congruences would follow if we could prove that
\begin{equation}
\label{eq:ppc}
pp(n)\equiv 0\ (\mathrm{mod}\ 5)\ \mathrm{if}\ n\equiv 2,3,4\ (\mathrm{mod}\ 5).
\end{equation}
This is because of the result
\[p_j(n)=pp\left(\frac{n-j(2j-1)}{2}\right)\]
of Theorem 1 above. There are 15 cases to consider, but fortunately they can all be done at once.

We want to prove that for each of the above pairs $(n,j)$ mod 5, the quantity $(n-j(2j-1))/2$ is either not an integer or else is 2, 3 or 4 mod 5. For it to be an integer we must have $j\equiv n$ mod 2. Hence we have a pair $(n,j)$ which modulo 5 have given values $(n',j')$, say, and are such that $j\equiv n$ mod 2. This means that
\[n=5s+5j'-4n'+10t, \ \mathrm{and}\ j=5s+j',\]
for some integers $s,t$. But then
\begin{equation}\label{eq:cnd}\frac{n-j(2j-1)}{2}\equiv 3j'-2n'-j'^2\ (\mathrm{mod}\,5).\end{equation}
Thus, to prove that (\ref{eq:ppc}) imply all of (\ref{eq:p0})--(\ref{eq:p4}) we need only verify that for each of the 15 pairs $(n',j')$
\[(0,1),(0,2),(1,0),(1,3),(1,4),(2,1),(2,2),(2,4),(3,0),(3,3),(4,\mathrm{all}),\]
mod 5 it is true that the right side of (\ref{eq:cnd}) is 2, 3 or 4 mod 5, which is a trivial exercise. $\Box$

It remains to establish (\ref{eq:ppc}). We have, modulo 5,
\[\frac{1}{(1-t)^2}\equiv \frac{(1-t)^3}{(1-t^5)},\]
and therefore
\[\prod_{j\ge 1}\frac{1}{(1-x^j)^2}\equiv \frac{\prod_{j\ge 1}(1-x^j)^3}{\prod_{j\ge 1}(1-x^{5j})}.\]
On the other hand it is known that
\[\prod_{j\ge 1}(1-x^j)^3=\sum_{n\ge 0}(-1)^n(2n+1)x^{{n+1\choose 2}}.\]
Consequently,
\[\sum_{k\ge 0}pp(k)x^k\equiv \left(\sum_{n\ge 0}(-1)^n(2n+1)x^{{n+1\choose 2}}\right)\left(\sum_{m\ge 0}p(m)x^{5m}\right).\]
Now all exponents of $x$ on the right are of the form $5m+{n+1\choose 2}$. Since ${n+1\choose 2}$ is always 0,1, or 3 mod 5, we have surely that $pp(k)\equiv 0$ if $k\equiv 2,4$ mod 5. Finally, if ${n+1\choose 2}\equiv 3$ mod 5, then $n\equiv 2$, so $2n+1\equiv 0$, and again the coefficient of $x^k$ vanishes mod 5. $\Box$

\vspace{.5in}

{\noindent\small Center for Combinatorics, LPMC, Nankai University, Tianjin 300071, P. R. China\\
\texttt{<chen@nankai.edu.cn>}}

\smallskip

{\noindent\small Center for Combinatorics, LPMC, Nankai University, Tianjin 300071, P. R. China\\
\texttt{<ji@nankai.edu.cn>}}

\smallskip

{\noindent\small Department of Mathematics, University of
Pennsylvania,
Philadelphia, PA 19104, USA\\
\texttt{<wilf@math.upenn.edu>}}

\end{document}